\documentstyle[12pt]{article}
 \begin{document}
 \baselineskip=15pt
 \newcommand{\ekq}{e^{2\pi i\frac{k}{q}}}
 \newcommand{\rf}{Ramanujan - Fourier~}
 \newcommand{\ekqn}{e^{2\pi i\frac{k}{q}n}}
 \newcommand{\ekqm}{e^{-2\pi i\frac{k}{q}}}
 \newcommand{\ekqmn}{e^{-2\pi i\frac{k}{q}n}}
 \newcommand{\wk}{Wiener - Khintchine formula~}
 \newcommand{\ekqp}{e^{2\pi i\frac{k'}{q'}}}
 \newcommand{\di}{\displaystyle}
\newcommand{\NI }{\noindent}
\newcommand{\be}{\begin{equation}}
\newcommand{\ee}{\end{equation}}
\begin{center}
 {\Large \bf LINKING THE CIRCLE AND THE SIEVE: 
RAMANUJAN- FOURIER SERIES}\\
 \vspace{1cm}
 {H. Gopalkrishna Gadiyar and R. Padma\\
AU-KBC Research Centre, M.I.T. Campus of Anna University\\
Chromepet, Chennai 600 044 India\\
E-mail: \{gadiyar, padma\}@au-kbc.org}
\end{center}

 \vspace{2cm}
 \begin{center}
 {\bf Abstract}
 \end{center}

Currently the circle and the sieve methods are the key tools in analytic number theory. In this paper the unifying theme of the two methods is shown to be \rf series.

\newpage
\NI  {\bf Introduction.} The two well known methods in additive number theory are the circle method and the sieve method. The circle method is based on using a generating function (See Section 3) and noting along with Ramanujan and Hardy that the rational points on the circle contribute most and then through estimates showing that the contribution from the other points is small. The other method is the sieve method which looks entirely different and tries to use probabilistic arguments by roughly giving a density of occurrence to multiples of numbers. In this note we show that the \rf series unite two different streams of thought, the analytic and the probabilistic. 

To explain this in an intuitive way, we consider the twin prime and the related Goldbach problem and break the paper into three parts: psychological, logical and chronological. The psychological part conveys the intuition behind the argument. The logical part displays the techniques needed. The chronological part gives a whiggish(!) history of the subject (as far as we could trace with available resources).  

\NI  {\bf 1. Psychological.} Let us assume that the primes follow some distribution. Then the twin prime problem would be deducible by studying the autocorrelation function of that distribution. The argument would require that the distribution have a Fourier series. 

An autocorrelation function of a process is defined to be the expectation of the product of two random variables obtained by observing the process at different times. The well-known \wk states a relationship between two important characteristics of a random process: the power spectrum of the process and the correlation function of the process [12]. It is used practically to extract hidden periodicities in seemingly random phenomena. The twin prime problem asks whether there are infinitely many prime pairs of the form $(p,~p+h)$ and how these pairs are distributed. One immediately notes that this is a problem of finding autocorrelation of primes and hence there should be a corresponding \wk .

\NI {\bf 2. Logical.} We state here the conventional \wk , explain \rf series of an arithmetical function and the \wk for an arithmetical function which leads to the conjecture of Hardy and Littlewood on the twin prime problem.  We also give the stochastic interpretation of the Riemann zeta-function due to Rota.

\NI {\bf 2.1~The Wiener - Khntchine formula. } The \wk [12] basically says that if
 \be
 f(t) ~=~ \sum_{n} f_n e^{i\lambda_n t} \, ,
 \ee
 then
 \be
 \lim_{T \rightarrow \infty} \frac{1}{2T} \int_{-T}^T f(t+\tau)\overline{f(t)} dt
 = \sum_{n} \vert f_n \vert ^2 e^{i \lambda _n \tau} \, . 
 \ee
 The left hand side of (2) is called an autocorrelation function. The right hand side is nothing but the power spectrum. 

\NI {\bf 2.2~\rf series.} The \rf series of an arithmetical function $a(n)$ is an expansion of the form
 \be
 a(n) ~=~ \sum_{q=1}^{\infty} a_q c_q(n) \, ,
 \ee
 where
 \be
 c_q(n) ~=~ \sum_{\stackrel{k=1}{(k,q)=1}}^q \ekqn \, ,
 \ee
 and $(k,q)$ denotes the greatest common divisor of $k$ and $q$. $c_q(n)$ is known as the Ramanujan sum and $a_q$'s are called the \rf coefficients. The \rf coefficient $a_q$'s are evaluated as follows. Denote by $M(g)$ the mean value of an arithmetical function $g$, that is,
 \be
 M(g) ~=~ \lim_{N \rightarrow \infty} \frac{1}{N} \sum_{n \le N} g(n) \, .
 \ee
 For $1 \le k \le q, ~ (k,q)=1,$ let
 $e_{\frac{k}{q}}(n)=e^{2\pi i \frac{k}{q} n}$, ($n \in \cal N$). If
 $a(n)$ is an arithmetical function with expansion (3),
 then
 \be
 a_q ~=~ \frac{1}{\phi (q)} M(a ~c_q) ~= \frac{1}{\phi (q)} ~\lim_{N \rightarrow
 \infty} \frac{1}{N}
 \sum_{n \le N} a(n) c_q(n)   ~,
 \ee
where $\phi (q)$ denotes the Euler totient function. Also,
 \be
 M(e_{\frac{k}{q}}  ~\overline{e_{\frac{k'}{q'}}}) = \left \{ \begin{array}{ll}
 1 \, , &if~~ \frac{k}{q}~=~ \frac{k'}{q'} \, ,\\
 0 \, , & if~~ \frac{k}{q} ~\ne ~\frac{k'}{q'} \, . 
 \end{array} \right .
 \ee
 The Ramanujan sum $c_q(n)$ satisfies some nice arithmetical properties which are stated below.

 \noindent (i) $c_q(n)$ is multiplicative. That is,
 \be
 c_{qq'}(n)~=~ c_q(n)~c_{q'}(n)\, {\rm ~if~} (q,q')=1 \, . 
 \ee

 \noindent (ii) If $p$ is prime, then
 \be
 c_p(n)~=~ \left \{ \begin{array}{ll} -1 \, , &{\rm ~if ~} p \not | n \, ,
 \\
 p-1, & {\rm ~if~} p | n \, ,
 \end{array} \right . 
 \ee
 where $a | b$ means $a$ divides $b$ and $a \not | b$ means $a$
 does not divide $b$.  

In number theory, it is  customary to use the von Mangoldt function $\Lambda
 (n)$ in the place of characteristic function of primes where $\Lambda
 (n)$ is defined as follows.
\be
\Lambda (n)=\left \{ \begin{array}{ll}
 \log p \, , &{\rm if}~ n~=~ p^k\, ,~ p, ~{\rm a~prime}~k,~{\rm a~positive~integer}\\
 0 \, , & ~{\rm otherwise}
 \end{array} \right .
\ee
 The \rf expansion for ${\di \frac{\phi(n)}{n}} \Lambda
 (n)$ can be obtained using the properties of $c_q(n)$ and it is given by 
\be
 \frac{\phi(n)}{n} \Lambda (n) = \sum_{q=1}^\infty \frac{\mu(q)}{\phi(q)}
 c_q(n) \, , 
\ee
where $\mu (q)$ is the M\"{o}bius function. 
 
\NI Also, for a given integer $h$,
 \be
 C(h) \stackrel{\rm def}{=} \sum_{q=1}^\infty \frac{\mu ^2 (q)}{\phi ^2(q)}
 c_q(h)=
 \left \{ \begin{array}{ll}
 2 {\di \prod_{p>2} \left (1-\frac{1}{(p-1)^2}\right )} ~
 {\di \prod_{\stackrel{p | h}{p >2}} \left (\frac{p-1}{p-2}\right )} \, ,
 & {\rm ~if~} h {\rm ~is~ even,}\\
 0 \, , & {\rm ~if~}h {\rm~is~odd,}
 \end{array} \right . 
 \ee
 where ${\di \prod_p}$ denotes the product over primes.

 The proof of (12) follows from the multiplicative properties of
 $\mu (q),~\phi (q)$ and $c_q$. For a given $h$, the series on the
 left hand side of (12) is absolutely convergent and hence has the
 Euler product expansion
 \be
 \prod_{p} \left ( 1+\frac{\mu ^2 (p)}{\phi ^2 (p)} c_p(h) \right )
 = \prod_{p ~\mid\!\!/~h} \left ( 1-\frac{1}{(p-1)^2} \right )
 \prod_{p | h} \left ( 1+\frac{1}{p-1} \right ) \, , 
 \ee
 where we have used the property (ii) of $c_q$ given by (9). When $h$ is odd, the infinite product on the right hand side of
 (13) is $0$
 and when $h$ is even, it is equal to
 \be
 2~\prod_{p > 2} \left ( 1-\frac{1}{(p-1)^2} \right )
 \prod_{\stackrel{p | h}{p>2}} \left ( 1+\frac{1}{p-1} \right )
 \left (1-\frac{1}{(p-1)^2} \right )^{-1} \, ,
 \ee
 which on simplification gives the right hand side of (12).

\NI {\bf 2.3 The \wk for an arithmetical function.} The \wk for an arithmetical function $a(n)$ having the \rf series (3)  can be stated as follows.
 \be
 \lim_{N \rightarrow \infty} \frac{1}{N} \sum_{n \le N} a(n)a(n+h) =
 \sum_{q=1}^{\infty} a_q^2 c_q(h) \, .
 \ee
(15) holds for those arithmetical functions whose \rf series are  absolutely and uniformly convergent. But the \rf series for
 ${\di \frac{\phi (n)}{n}} \Lambda (n)$
 does not belong to this class of functions.
 However, if we assume the truth of the theorem in the most general case, it would give the formula conjectured by Hardy and Littlewood [9] which we now state.

 {\it There are infinitely many prime pairs $p,
 p+h$ for every even integer $h$ and if $\pi _h(N)$ denotes the number of prime
 pairs
 less than $N$, then
 \be
 \pi _h(N) \sim C(h)
 \frac{N}{\log^2N} \, . 
 \ee
 }
where $C(h)$ is given by (12).

The \wk for ${\di \frac{\phi (n)}{n}} \Lambda (n)$ is given by
 \begin{eqnarray}
 \lim_{N \rightarrow \infty} \frac{1}{N}\sum_{n \le N} \frac{\phi (n)}{n} \Lambda
 (n)
 \frac{\phi (n+h)}{n+h} \Lambda (n+h)& =&
 \sum_{q=1}^\infty \frac{\mu ^2 (q)}{\phi ^2(q)} c_q(h)\\
 &=& C(h) \, . 
 \end{eqnarray}
 Let
 \be
 \Psi (h,N) =
 \sum_{n \le N} \frac{\phi (n)}{n} \Lambda (n)
 \frac{\phi (n+h)}{n+h} \Lambda (n+h) \, .
 \ee
 Then (17) implies that
 \be
 {\di \frac{\Psi( h,N)}{N}} \sim C(h)  \, . 
 \ee
 The terms of the sum on the left hand side of (17) are non zero if and only if both $\Lambda (n)$ and $\Lambda (n+h)$
 are prime powers, say $p^k$ and $q^l$ respectively, for primes $p$ and $q$. Passing from (20)to (16) is a standard exercise in number theory. See [6]. 

We now give the strong numerical evidence based on computerized calculations which indicate the plausibility of the \wk [6].
 For $h$ = 2, 4 and 6, the formula (20) gives
 \begin{eqnarray*}
 \Psi(2,N) &\sim& C(2)~N\\
 \Psi(4,N) &\sim& C(4)~N\\
 \Psi(6,N) &\sim& C(6)~N
 \end{eqnarray*}
 Since $C(4)=C(2)$ and $C(6)=2C(2)$, there should be approximately equal numbers of prime power pairs differing by 2 and by 4, but about twice as many differing by 6. We have given below the actual values of $\Psi (h,N)$ and also the ratio
 $C(h) / \left ({\di \frac{\Psi(h,N)}{N}}\right )$ (fourth column) for $h~=~2,~4$ and 6 and $N$ upto $10^6$. We have used the value of
 $C(2) = 2~{\di \prod_{p > 2} \left ( 1-\frac{1}{(p-1)^2} \right )} \sim
 1.320323632$ to compute the ratio.

 \begin{center}
 {\bf Table 1}~~~~~~~~~~~~~~~
 \end{center}

 \begin{tabular}{|c|c|c|c|}
 \hline
 \hline
 N& $\Psi(2,N)$ & ${\di \frac{\Psi(2,N)}{N}} $ & Ratio\\
 \hline
 100000 & 131522.552204 & 1.315226 & 1.003876 \\
 200000 & 264287.347531 & 1.321437 & 0.999158 \\
 300000 & 393317.025988 & 1.311057 & 1.007068 \\
 400000 & 525523.270611 & 1.313808 & 1.004959 \\
 500000 & 654557.716460 & 1.309115 & 1.008562 \\
 600000 & 789035.163302 & 1.315059 &1.004004 \\
 700000 & 919941.157912& 1.314202 & 1.004658 \\
 800000 & 1049182.174335 & 1.311478& 1.006745 \\
 900000 & 1180813.946552 &1.312015 & 1.006332 \\
 1000000 & 1312843.985016 & 1.312844 & 1.005697\\
 \hline
 \end{tabular}

 \begin{center}
 {\bf Table 2}~~~~~~~~~~~~~~~
 \end{center}

 \begin{tabular}{|c|c|c|c|}
 \hline
 \hline
 N& $\Psi(4,N)$ & ${\di \frac{\Psi(4,N)}{N}} $ & Ratio\\
 \hline
 100000 & 130212.335085 & 1.302123 & 1.013977 \\
 200000 & 260492.247225 & 1.302461 & 1.013714 \\
 300000 & 390320.617781 & 1.301069 & 1.014799 \\
 400000 & 527155.226011 & 1.317888 & 1.001848 \\
 500000 & 653649.051733 & 1.307298 & 1.009964 \\
 600000 & 789177.513123& 1.315296 &1.003823 \\
 700000 & 923982.224287 & 1.319975 & 1.000264 \\
 800000 & 1054670.388142 & 1.318338& 1.001506 \\
 900000 & 1180133.117590 & 1.311259 & 1.006913 \\
 1000000 & 1307978.775955 & 1.307979 & 1.009438\\
 \hline
 \end{tabular}

  \begin{center}
 {\bf Table 3}~~~~~~~~~~~~~~~
 \end{center}

 \begin{tabular}{|c|c|c|c|}
 \hline
 \hline
 N& $\Psi(6,N)$ & ${\di \frac{\Psi(6,N)}{N}} $ & Ratio\\
 \hline
 100000 & 261289.742091 & 2.612897 & 1.010620 \\
 200000 & 523391.109218 & 2.616956 & 1.009053 \\
 300000 & 787393.641752 & 2.624645 & 1.006097 \\
 400000 & 1056087.319082 & 2.640218& 1.000162 \\
 500000 & 1316336.875799 & 2.632674 & 1.003029 \\
 600000 & 1579274.310330 & 2.632124 &1.003238 \\
 700000 & 1839327.388416 & 2.627611 & 1.004961 \\
 800000 & 2104826.034045 & 2.631033 & 1.003654 \\
 900000 & 2368450.398104 & 2.631612 & 1.003434 \\
 1000000 & 2631198.406265 & 2.631198 & 1.003591\\
 \hline
 \end{tabular}
\vspace{0.5cm}

\NI {\bf 2.4 Stochastic interpretation of Rota.} In this section we point out the fact that Rota's stochastic interpretation of the Riemann - zeta function [19], [1] is based on profinite characters. We will see that Rota considers the group $C_{\infty}$ of rational numbers modulo 1 and crucially bases his arguments (given below) on $C_{\infty}^*$, the group of characters of $C_{\infty}$. {\it Both Ramanujan and Rota are using the same tool as $e^{2\pi i \frac{k}{q}n}$ of Ramanujan is a concrete realization of the characters of the profinite group of Rota}. It seems that Rota did not notice the Ramanujan connection. As Sieve methods [7] are based on probabilistic considerations, the circle and the sieve methods are linked through the ideas of Rota and Ramanujan.

Consider $A$ a subset of positive integers $\cal N $, then the arithmetic density is defined as 
\be
dens(A) = \lim_{n \rightarrow \infty} \frac{1}{n} || A \cap \{ 1,2,...,n\} || \, ,
\ee 
whenever the limit exists. It is immediately obvious that the dens($\cal N $)=1 and $dens(A_p)=\frac{1}{p}$ where $A_p$ is the set of multiples of $p$. For technical reasons (lack of countable additivity) the right way to go about it is to choose a number $s > 1$ define the measure of a positive integer $n$ to be $\frac{1}{n^s}$. Then it turns out that the measure of $\cal N $ is equal to $\di{\zeta (s) = \sum_{n=1}^\infty \frac{1}{n^s}}$. Hence a countably additive measure $P_s$ on the set $\cal N$ defined as
\be
P_s(A) = \frac{1}{\zeta (s)} \sum_{n \in A} \frac{1}{n^s} \, .
\ee
Further the fundamental property 
\be
P_s(A_p \cap A_q) = P_s(A_p) P_s(A_q) = \frac{1}{p^sq^s} 
\ee
can be checked and ${\di \lim_{s\rightarrow 1} P_s(A) = dens(A)}.$
That is, the arithmetic density though not a probability is the limit of probabilities. 

Rota then gives a combinatorial twist to the problem. He considers a cyclic group $C_r$ of order $r$. Every character $\chi$ of the group $C_r$ has a kernel which is a subgroup of $C_r$. Further every sequence $\chi_1 ,..., \chi_s$ of characters of $C_r$ has a joint kernel which is a subgroup of $C_r$. By a joint kernel of a sequence of characters we mean the intersection of their kernels. Suppose one were to choose a sequence of $s$ characters independently and randomly and ask the question: what is the probability that the joint kernel is a certain subgroup $C_n$ of $C_r$? It can be seen that with probability $\frac{1}{n}$ the kernel of a randomly chosen character will contain the subgroup $C_n$ as there are $r$ characters of the group $C_r$ and $\frac{r}{n}$ such characters will vanish on $C_n$. So the probability of the joint kernel will contain $C_n$ is $\frac{1}{n^s}$. Let $P_{C_n}$ denote the probability that the joint kernel of characters shall be $C_n$. It follows that ${\di
\frac{1}{n^s} = \sum_{n |d|r} P_{C_d}}$
This is based on the fact that the partially ordered set of subgroups of a cyclic group $C_r$ and the partially ordered set of the divisors of r are isomorphic. Next using M\"{o}bius inversion and change of variable  $d= nj$ it follows that 
\be
P_{C_n} = \sum_{n|d|r} \mu(\frac{d}{n}) \frac{1}{d^s}
\ee
\be
P_{C_n} = \frac{1}{n^s} \sum_j \mu(j) \frac{1}{j^s}
\ee
The variable $j$ ranges over a subset of divisors of $r$. If the sum ranged over all positive integers $j$ we would get the probabilistic interpretation 
\be
\frac{1}{n^s} \frac{1}{\zeta(s)} \, .
\ee
This is done by replacing the finite cyclic group $C_n$ by a profinite cyclic group. Take the group $C_{\infty}$ of rational numbers modulo 1. For every positive integer $n$ the group $C_{\infty}$ has unique finite subgroup $C_n$ of order $n$. The character group $C_{\infty}^*$ of $C_{\infty}$ is a compact group, and the Haar measure is a probability measure. The group $C_{\infty}^*$ is the desired profinite group. Mimicking the argument made earlier it can be seen that
\be
 \frac{1}{n^s}=\sum_{n |d} P_{C_d}
\ee
and by M\"{o}bius inversion
\be
P_{C_n} = \sum_{n|d} \mu(\frac{d}{n}) \frac{1}{d^s} = \frac{1}{n^s} \frac{1}{\zeta(s)} \, .
\ee 
The characters of $C_{\infty}$ are given by $e^{2\pi i\frac{k}{q}n}$ and so one immediately observes that the connection between $C_{\infty}^*$ of Rota and the \rf expansion.

\NI  {\bf 3. Chronological.}  In this section we will re-analyze the circle method and extract from it the essence which leads to great simplification and an intuitive understanding of
 the issues involved. We will also give the historical roots of probabilistic ideas in Hardy's work. The two streams will merge because of Rota's profinite characters being related to \rf series. 

\NI {\bf 3.1 Meaning of the circle method.}  We begin by retracing what we feel should have been the train of thought
 of Hardy, Littlewood and Ramanujan, as far as it is possible, from their
 published researches. The first problem to be attacked by the circle method
 was the partition problem. As is well known, this reduces to understanding
 the generating function
 \be
 \sum_{n=0}^{\infty} p(n)\,x^n ~=~ \frac{1}{(1-x)(1-x^2)(1-x^3)......} \, .
 \ee
The key observation that led to the development of the circle method was
 that $p(n)$ could be written as
 \be
 \frac{1}{2\pi i} \int_C{\frac{1}{(1-z)(1-z^2)(1-z^3).....} z^{-n-1}dz} \, ,
\ee 
where $C$ is the circle $|z|=r$ and $r<1$. Hardy and Ramanujan observed that the unit circle is covered by infinity of singularities corresponding to the poles of the generating function at all the rational points on the unit circle. This led them to carry out some extremely subtle and delicate analysis which gave remarkably accurate numerical agreement with available results. Having got this success, Hardy along with Littlewood went on to attack the other problems in additive number theory like the Waring's problem and the Goldbach problem. In each case they got approximate formulae in terms of what is called the singular series.
 At this point, Hardy and Littlewood seem to have shifted their focus to the
 complex analytic aspects of the circle method whereas Ramanujan wrote what is
 probably the one of the deepest papers in number theory in which he introduced
 the concept of what is now called the Ramanujan - Fourier expansion [18].

 What Ramanujan does is to show by simple, yet ingenious methods that a wide
 range of arithmetical functions have Ramanujan - Fourier expansions, that is, an expansion of the form (3). He however does not indicate any formula for getting the Ramanujan -
 Fourier coefficients $a_q$ which are the backbone of Fourier analysis.
 This was done by R. D. Carmichael [3] a little later. He showed that
 the \rf coefficients $a_q$ can be got by the formula (6).
  M. Kac, E. R. Van Kampen and A. Wintner [10], [11] have pointed out the almost periodic nature of the Ramanujan-Fourier expansions. See also [13] and [20].

In the mean time, Hardy and Littlewood realised that their singular series could be evaluated because it was actually a Ramanujan - Fourier expansion. To do this, Hardy used the techniques he had developed in one of his most important papers [8] in which he stressed the multiplicative nature of $c_q(n)$ ((8) and (9)) which enabled him to get closed form formula for the singular series. Using the properties of $c_q$ Hardy also obtained a \rf expansion of ${\di \frac{\phi (n)}{n}} \Lambda (n)$ given by (11).

As historical aside, we would like to remark that clues of Ramanujan - Fourier expansion playing an important role in additive number theoretic problems are available in the work of Glaisher [9] in which a simple minded
 partial fraction expansion yields this.
 \begin{eqnarray}
 \frac{1}{(1-x)(1-x^2)(1-x^3)} &=& \frac{1}{6(1-x)^3} ~+~ \frac{1}{4(1-x)^2}\\
  &+& \frac{17}{72(1-x)} ~+~ \frac{1}{8(1+x)}\\
 &+& \frac{1}{9(1-\omega x)} ~+~ \frac{1}{9(1-\omega ^2x)} \, ,
 \end{eqnarray}
 where $\omega $ and $\omega ^2$ denote the two complex cube roots of unity.
 If
 \be
 \frac{1}{(1-x)(1-x^2)(1-x^3)}~=~ 1 ~+~ \sum_{n=1}^\infty r(n)x^n \, ,
\ee
 then
 \be
 r(n) ~=~ \frac{(n+3)^2}{12} -\frac{7}{72} + \frac{(-1)^n}{8} + \frac{2}{9} \cos
 (\frac{2n\pi}{3}).
 \ee

 Over the last half century, several refinements have been made in attacking additive number theory problems. In the circle method, as applied by Hardy and Littlewood to the Goldbach problem, strong hypotheses regarding primes in arithmetic progressions have to be made. However, ultimately what is got is a main term in the formula which is given in terms of the Ramanujan - Fourier series and an error term. The hard part is to control the error term. At this point, we would like to indicate what we consider the line that Ramanujan would have taken had he not lost interest in the Ramanujan - Fourier expansion.

H. Wilf [21] remarks that ``A generating function is a clothesline on which we hang up a sequence of numbers for display". But is a clothesline necessary, that is, do we need a generating function at all? Rather than considering the generating function and then applying the circle method which yields a Ramanujan - Fourier expansion, can we not directly get the answer? What should be done is a clean shave with Occam's Razor.

Consider a typical additive number theory problem. We would begin with
 a subset {\cal A} of positive integers and let $a(n)$ denote its characteristic function. If the problem is to find out $r(n)$ which is the number of ways $n$ can be written
 as a sum of $k$ elements of the set {\cal A}, we would extract the $n^{th}$  coefficient of $(\sum a(n) x^n)^k$. This corresponds to taking the $k$-fold convolution. In other words,
 $$
 r(n) ~=~ \sum_{i_1+i_2+...+i_k=n} a(i_1)a(i_2)...a(i_k).
 $$
 All the analysis carried out over the past few years finally leads to an approximate formula for this convolution.

 If $a(n)$ has the \rf expansion (3), then
 \be
  r(n) \sim  \left (\begin{array}{c} n+k-1 \\ k-1 \end{array}\right )
 \displaystyle {\sum_{q=1}^\infty a_q^k c_q(n) \, .}
 \ee
This formula is immediately understandable as the usual relationship
 between the convolution and multiplication which exists for Fourier series.
 The combinatorial factor which appears is due to the fact that for almost periodic functions [2] only approximate formulae exist and this combinatorial factor is necessary there. 

 The reason why the rational points dominate is simply because
 of the fact that the $a(n)$ have
 Ramanujan - Fourier expansions which give rise to the simple poles at
 all the rational points on the unit circle. For,
 \begin{eqnarray}
 \sum_{n=1}^\infty a(n)x^n &=& \sum_{n=1}^\infty
 \sum_{q=1}^\infty \sum_{\stackrel{k=1}{(k,q)=1}}^q a_q e^{2\pi i
 \frac{k}{q}n}x^n\\
 &=& \sum_{q=1}^\infty \sum_{\stackrel{k=1}{(k,q)=1}}^q a_q \sum_{n=1}^\infty
 \left(e^{2\pi i \frac{k}{q}}x\right)^n\\
 &=& \sum_{q=1}^\infty \sum_{\stackrel{k=1}{(k,q)=1}}^q a_q \frac{e^{2\pi i
 \frac{k}{q}}x}{1-e^{2\pi i\frac{k}{q}}x} \, .
 \end{eqnarray}
 
 For example, for the Goldbach problem, it is known that $\di{\frac{\phi (n) \Lambda(n)}{n}}$ has a \rf expansion (11).
  It is this $ \displaystyle{\frac{\mu (q)}{\phi (q)}}$which is got by using the properties of primes
 in arithmetic progressions in the circle method. Hence if an approximate
 formula for convolution can be proved rigorously, the Goldbach problem
 can be solved. In other methods there are many sources of error and the
 methods do not tell us why things work. Further they involve many hypotheses. Hence a return to the ideas of Ramanujan seems to be absolutely necessary to clarify and simplify the circle method.

Let us now look at the Goldbach conjecture more closely. If $r(n)$ denotes the number of ways $n=2m$ can be expressed  as a sum of two primes, then from (11) and (36), 
\be
r(n) \sim C(n) \frac{n}{\log ^2n} \, ,
\ee
where $C(n)$ is given by (12). This formula was conjectured by Hardy and Littlewood. But Sylvester was the first mathematician to suggest an asymptotic formula for $r(n)$. He conjectured that
\be
r(n) \sim 2\pi (n) \prod \left (\frac{p-2}{p-1}\right )\, ,
\ee
where $\pi (n)$ denotes the number of primes up to n and the product extends over all odd primes $p\le n$ and $p \not |n$. There is no evidence in the literature as to how he arrived at this formula. We know by prime number theorem that
\be
\pi (n) \sim \frac{n}{\log n}\, ,
\ee
and by Merten's formula
\be
\prod_{p \le x} \left ( 1-\frac{1}{p} \right ) \sim \frac{e^{-\gamma}}{\log x} \, ,
\ee
where $\gamma $ is Euler's constant. This shows that
\be
\prod_{p\le \sqrt{n}} \left (\frac{p-2}{p-1}\right ) = \prod_{p\le \sqrt{n}} \left (1-\frac{1}{(p-1)^2}\right )\prod_{p \le \sqrt{n}} \left ( 1-\frac{1}{p} \right ) \sim \frac{2Ae^{-\gamma }}{\log n} \, ,
\ee
where 
\be
A=\prod_{p > 2} \left (1-\frac{1}{(p-1)^2}\right )\, ,
\ee
where the product runs over all odd primes. Thus
\be
r(n) \sim \frac{4Ae^{-\gamma }n}{\log ^2n}  \prod_{\stackrel{p | n}{p >2}}\left (\frac{p-2}{p-1}\right )\, .
\ee
One may call this as `Sylvester's formula'. This formula is clearly wrong. It is right in its most important terms, but the constant $4Ae^{-\gamma }$ is not correct and the correct formula should just have $2A$ as the constant factor. However, the sieve argument used by Merlin and Brun gives an asymptotic formula similar to the one in given by (46) which is given below.

\NI {\bf 3.2 Probabilistic Interpretation of the Goldbach conjecture.} Here one forms the table
$$
\begin{array}{cccccccccc}
1,&2,&3,&,4,&5,&6,&7,&\cdots ,&n-1&~(m)\\
n-1,&n-2,&n-3,&,n-4,&n-5,&n-6,&n-7,&\cdots ,&1&~(m') 
\end{array}
$$
where $n$ is an even number and the table read vertically gives the decompositions $n=m+m'$ of $n$ into positive integers $m$ and $m'$. Let us now perform the sieve of Eratosthenes on both rows of the table sieving with the first $l$ primes starting from the right of the lower row. A decomposition $m+m'$ is considered erased when either $m$ or $m'$ is erased.

Two cases occur here. If the prime $p|n$, then the erasures in the second row fall immediately below the corresponding erasures in the first, and the corresponding erasures in the first, and the number of decompositions erased is approximately $\frac{n}{p}$. If $p \not | n$, then the erasures never correspond, and now the number of decompositions is approximately  $\frac{2n}{p}$. If $r_l(n)$ corresponds to the number of decompositions $n=m+m'$ such that $m$ or $m'$ are not divisible by any of the first $l$ primes, we have
\be
r_l(n)=n \prod_{p|n,~p\le p_l} \left (1-\frac{1}{p}\right )
\prod_{p\not | n,~p\le p_l} \left (1-\frac{2}{p}\right )\, .
\ee
This formula is correct for a fixed $l$. But if we assume the truth of this formula for $p_l\sim \sqrt{n}$, then we have $r(n)=r_l(n)$ and
\be
r(n)=\frac{1}{2} n \prod_{p|n} \left (\frac{p-1}{p-2}\right )
\prod_{p\le \sqrt{n}} \left (1-\frac{2}{p}\right )\, ,
\ee
where $p=2$ is now excluded. By Merten's formula, if $3 \le p \le \sqrt{n}$, 
\be
\prod \left (1-\frac{2}{p}\right )= \prod \left (1-\frac{1}{(p-1)^2}\right )\left (1-\frac{1}{p}\right )^2 \sim \frac{16Ae^{-2\gamma }}{\log ^2n} \, .
\ee
Thus,
\be
r(n) \sim  \frac{8Ae^{-2\gamma }n}{\log ^2n}\prod_{p|n}\left (\frac{p-1}{p-2}\right )\, .
\ee
This is the formula to which Brun's argument naturally leads. Like Sylvester's formula this is also wrong but correct up to the factor of $4e^{-2\gamma }$. 

A similar probabilistic argument for the twin prime conjecture has been lucidly explained by Polya [17].

Let us now recall the probabilistic interpretation of $\zeta (s)$ given in Section 2.3 and note that the $e^{2\pi i \frac{k}{q}n}$ of Ramanujan is a concrete realization of the characters of the profinite group of Rota. Thus the ideas of circle method and sieve method can be linked through the work of Rota via \rf expansion.

\NI {\bf 3.3 Related ideas.} Physicists have carried out statistical analysis of primes recently and shown that the distribution of primes if related to $1/f$ noise. Marek Wolf [22] has shown that the power spectrum displays the $1/ƒ^{\beta }$ behaviour with the exponent $\beta \sim 1.64 $ by performing the Fourier transform of the "signal" given by the number of primes contained in the successive intervals of equal length $l = 2^{16} = 65536$ up to $N = 2^{38} \sim  2.749 × 10^{11}$. Here this slope $\beta $ does not depend on the length of the sampled intervals, which suggests some kind of self-similarity in the distribution of primes. From Hardy and Ramanujan we know that primes have \rf expansion. Michel Planat [16] has used \rf series to signal processing (which is unrelated to number  theory.)

Before we conclude we make two comments of historical nature.
\begin{enumerate}
\item In a paper on time-series analysis [5], A. Einstein has discussed the autocorrelation function and its relationship to the spectral content of a time series many years before Wiener and Khintchine.
\item It is well known that Hardy had true antipathy to probability theory. In [4] Diaconis shows that despite his antipathy, Hardy contributed significantly to modern probability. ``His work with Ramanujan begat probabilistic number theory. His work on Tauberian theorems
and divergent series has probabilistic proofs and interpretations. Finally Hardy spaces are a
central ingredient in stochastic calculus."
He refers to Hardy's work with Littlewood on the Goldbach and prime k-tuples conjecture as an instance of his insightful probabilistic thinking. ``Hardy's papers give a sophisticated development of conjectured asymptotics which suggest that at least one of the authors was quite familiar with probabilistic heuristics. .... Hardy refers to what I would call 'probabilistic reasoning' as
"a priori judgment of common sense" in his expository account (Hardy (1922, page 2)."
\end{enumerate}
Thus the two streams of thought of Hardy one, the heuristic probabilistic (sieve) versus the other, the rigorous analytic (circle) are merged into one stream via \rf series.

\NI {\bf Conclusion: Ontological and epistemological.} The reader might have noticed that though there is no proof in the strict mathematical sense of the word, there is convincing calculation which shows the plausibility of the conjecture. The philosophy of this paper is what Roddam Narasimha [15] calls computational positivism. Here we have a model for the primes s behaving like random variables and based on this model a numerical calculation is carried out which agrees with ``experimental" data. This is the typical method used repeatedly in physics where one guesses a model and carries out a numerical calculation.

Both the circle method and the sieve method lead to \rf series. In the circle method, poles occurring at rational points lead to \rf series. In the sieve method, the profinite group characters give the right density on the integers. These two methods are therefore linked through \rf series.

\NI {\bf Dedication.} It is a pleasure to write an article in honour of Professor K. Ramachandra on his $70^{th}$ birthday. He belongs to the pre-independence nationalist school of scientists like Ramanujan (whom he admires.) He has been a symbol of resistance to neo-colonial fashions and an advocate of indigenous styles of research. 

\noindent {\bf References}
\begin{enumerate}
\item  K. S. Alexander, K. Baclawski, G.C. Rota, {\it A stochastic interpretation of the Riemann zeta-function}, Proc. Natl. Acad. Sci. USA, {\bf 90}, (1993) 697-99. 
\item H. Bohr, {\it Almost periodic functions,} Chelsea Publishing
 Company, New York, 1951.
\item R. D. Carmichael,  Expansions of arithmetical functions
 in infinite series, Proc. London Math. Soc. (2) 34 (1932), 1-26.
\item P. Diaconis, {\it G.H. Hardy and Probability????} Bull London Math. Soc. {\bf 34}, (2002) 385-402.
\item A. Einstein, {\it Method for the determination of the statistical values of observations concerning quantities subject to irregular fluctuations}, Archives des Sciences et Naturelles, {\bf 37}, (1914) 254-256.
\item H. G. Gadiyar and R. Padma, {\it Ramanjuan-Fourier series, the Wiener-Khintchine formula and the distribution of prime pairs}, Physica {\bf A} 269(1999) 503-510.
\item G. Greaves, {\it Sieves in Number Theory}, Springer-Verlag, Berlin-Heidelberg, 2001.
\item G. H. Hardy, Note on Ramanujan's trigonometrical function
 $c_q(n)$ and certain series of arithmetical functions,  Proc. Camb. Phil. Soc. {\bf 20}, (1921) 263 - 271.
\item G. H. Hardy and J. E. Littlewood, `Some problems of Partition
 Numerorum'; III: On the expression of a number as a sum of primes,  Acta Math. {\bf 44}, (1922) 1-70.
\item M. Kac, {\it Almost periodicity and the representation of
 integers as sums of squares,} Am. J. Math. {\bf 62}, (1940) 122 - 126.
 \item  M. Kac, E. R. Van Kampen and A. Wintner, {\it Ramanujan sums
 and almost periodic functions,} Am. J. Math. {\bf 62}, (1940) 107 - 114.
 \item  C. Kittel, {\it Elementary Statistical Physics}, John Wiley \& Sons Inc., New York, 1958.
\item J. Knopfmacher, {\it Abstract Analytic Number Theory}, North-Holland/Americal Elsevier, 1975.
\item J. Kubilius, {\it Probabilistic Methods in the Theory of Numbers}, Vol. 11, Translations of Mathematical Monographs, American Mathematical Society Providence, Rhode Island, 1964.
\item Roddam Narasimha, {\it Axiomatism and computational positivism - Two mathematical cultures in pursuit of exact sciences}, Economic and Political Weekly, August 30, 2003, 3650-3656. 
\item M. Planat,  H. Rosu and S. Perrine, {\it Ramanujan sum for signal processing of low frequency noise} Phys. Rev. E {\bf 66} (2002) 56128.
\item G. Polya, {\it The Random Walks of George Polya}, Appendix 12, Heuristic Reasoning in the Theory of Numbers, 265-275. 
\item  S. Ramanujan, On certain trigonometrical sums and their
 applications in the theory of numbers, Trans. Camb. Phil. Soc.
 22 (1918), 259 - 276.
\item G.C. Rota, {\it Combinatorial snapshots}, The Mathematical Intelligencer, {\bf 21} (1999) No.2, 8-14.
\item W. Schwarz and J. Spilker, {\it Mean values of Ramanujan expansions of almost even arithmetical functions}, In: Proc. 1974, Colloq. on Number Theory, Debrecen.
\item H. Wilf, {\it Generatingfunctionology,} Academic Press,
 1990.
\item  M. Wolf, $1/f$ noise in the distribution of prime numbers,
 Physica {\bf A} 241 (1997), 493-499.
 \end{enumerate}
 \end{document}